\documentclass[a4paper]{article}

\usepackage[version=3]{mhchem}
\usepackage{authblk}

\usepackage{amsthm,amssymb}
\usepackage{graphicx,url}
\newcommand{\cc}{concentration}
\newcommand{\emp}{\textbf}
\newcommand{\mma}[1]{\textbf{\texttt{#1}}}
\newcommand{\mmares}[1]{\texttt{#1}}

\theoremstyle{definition}

\newtheorem{theorem}{Theorem}

\newtheorem{example}{Example}
\newtheorem{remark}{Remark}

\newtheorem{definition}{Definition}

\newcommand{\ikde}{induced kinetic differential equation}
\newcommand{\Vg}{Volpert graph}
\newcommand{\Vi}{Volpert index}
\newcommand{\FHJg}{Feinberg--Horn--Jackson graph}

\newcommand{\Reactants}{\mathrm{Reactants}}
\newcommand{\Products}{\mathrm{Products}}
\newcommand{\Producing}{\mathrm{Producing}}

\newcommand{\following}[1]{\textbf{next}(#1)}

\newcommand{\ReactionKinetics}{\mma{ReactionKinetics}}

\newcommand{\N}{\mathbb{N}}
\newcommand{\R}{\mathbb{R}}
\newcommand{\Rp}{\mathbb{R_+}}
\newcommand{\AAA}{\mathcal{A}}
\newcommand{\MM}{\mathcal{M}}
\newcommand{\RR}{\mathcal{R}}

\newcommand{\cb}{\mathbf{c}}
\newcommand{\yb}{\mathbf{y}}
\newcommand{\zb}{\mathbf{z}}
\newcommand{\alphab}{\boldsymbol{\alpha}}
\newcommand{\betab}{\boldsymbol{\beta}}

\newcommand{\T}{^{\top}}

\title{How to generate species with positive concentrations for all positive times?}

\author[1]{Seyed Mahdi Shavarani}
\author[$^{\ast}$,2,3]{J\'anos T\'oth} 
\author[4]{B\'ela Vizv\'ari}
\affil[1]{Alliance Manchester Business School, University of Manchester, Manchester, UK, United Kingdom}
\affil[2]{Department of Analysis, Budapest University of Technology and Economics, Budapest, Hungary.  Fax: +36 1 463 3172; Tel: +36 1 463 2314; E-mail: jtoth@math.bme.hu}
\affil[3]{Chemical Kinetics Laboratory, Institute of Chemistry, E\"otv\"os Lor\'and University, Budapest, Hungary}
\affil[4]{Department of Industrial Engineering,
Eastern Mediterranean University, Famagusta, North Cyprus, Turkey}

\begin{document}
\maketitle

\begin{abstract}
Given a reaction (network) we are looking for minimal sets of species starting from which all the species will have positive concentrations for all positive times in the domain of existence of the solution of the induced kinetic differential equation. We present three algorithms to solve the problem. The first one essentially checks all the possible subsets of the sets of species. This can obviously work for only a few dozen species because of combinatorial explosion. The second one is based on an integer programming reformulation of the problem. The third one walks around the state space of the problem in a systematic way and produces all the minimal sets of the advantageous initial species, and works also for large systems. All the algorithms rely heavily on the concept of Volpert indices, used earlier for the decomposition of overall reactions [Kov\'acs \textit{et al., Physical Chemistry Chemical Physics} 2004, \textbf{6}, 1236]. Relations to the permanence hypothesis, possible economic or medical uses of the solution of the problem are analyzed, and open problems are formulated at the end.
\end{abstract}

\section{Introduction}
Suppose we have a reaction (network), and our goal is that all the species have positive \cc s for all finite times in their domain of existence.
We would like to achieve this by putting only a few species (a minimal set of species) into the reaction vessel or test tube initially.
\begin{example}
Note that being positive even for all finite times does not mean that the stationary \cc\ of the given species will also be positive, as the simple example \ce{X -> Y} shows, where the stationary \cc\ of \ce{X} is zero, still it is positive for all finite times. This example also shows that a species with initially zero concentration can be positive for all positive times: if we start
form $x(0)>0, y(0)=0$ the concentration $y(t)$ will always be positive after time zero.
\end{example}
\begin{example}
The example \ce{2X -> 3X} shows that $x(0)>0$ implies positivity
for all positive times in the domain of solution, the interval $[0,1/x_0[$ which is bounded from above, i.e. \emp{not} the whole set of all nonnegative real numbers.
\end{example}

We shall introduce all our methods on a simple, still not uninteresting example.
\begin{example}
Let us consider the irreversible Michaelis--Menten reaction:
\begin{equation*}
  \ce{E + S <=> C -> E + P}
\end{equation*}
In this case it is intuitively obvious that all the species will have positive \cc s
if the initial \cc s of the following species are positive:
\begin{align*}
  &\{\ce{C}\},\\
  &\{\ce{E},\ce{S}\},\{\ce{E}, \ce{C}\},\{\ce{S}, \ce{C}\},\{\ce{C}, \ce{P}\},\\
  &\{\ce{E}, \ce{S}, \ce{C}\},\{\ce{E}, \ce{S}, \ce{P}\},\{\ce{E}, \ce{C}, \ce{P}\},\{\ce{S}, \ce{C}, \ce{P}\},\\
  &\{\ce{E}, \ce{S}, \ce{C}, \ce{P}\}.
\end{align*}
Of these possibilities the minimal ones are \{\ce{C}\} and \{\ce{E}, \ce{S}\}, in the sense that all other sets of species on the list are strict supersets of either set (or even both).
\end{example}
Rephrasing our goal is to create an algorithm producing minimal initial sets of species whose positive initial concentrations are enough to produce all the species in positive concentrations for all finite times. We are looking for effective algorithms
using as much chemical information as possible.

The structure of the paper is as follows.
In Section \ref{sec:problemform} we formally introduce the necessary concepts: we give the mathematical definition of reaction networks and their induced kinetic differential equation, the Volpert graph of the reaction, Volpert indexes of species and reaction steps. These definitions can be found in e.g.\ Ref. \cite{tothnagypapp}, but see also our previous paper\cite{kovacsvizvaririedeltoth}.

Section \ref{sec:bruteforce} contains a verbal description---together with further examples---of a refined version of the brute force solution.
Section \ref{sec:lp} provides arguments for an integer programming formulation.
To treat really large problems
we need a smart solution taking into consideration the structure of the problem.
This is given in Section \ref{sec:smart}.
Relevant parts of the codes for all the methods are presented either as pseudocodes, or as real codes.
More details can be found in our notebook written in the \mma{Wolfram language} or our Matlab worksheets available from the authors upon request.
Unusual it may be but not earlier than before the end shall we analyze possible applications and relations of the problem with other fields of kinetics; this would not have been easy to follow earlier. Note that our codes heavily use the \mma{Wolfram language} package \ReactionKinetics\ downloadable from the site \url{http://extras.springer.com}.

To end our introduction, we enumerate the notations used in the paper.

\begin{table}[!ht]
  \caption{Notations}
  \label{tbl:notations}
  \begin{tabular}{ll}
    \hline
    Notation & Meaning \\
    \hline
    $\AAA$&the set of atoms\\
    $\cb$&concentration vector at time $t$ \\
    $c_m(t)$ & concentration of species \ce{X($m$)} at time $t$ \\
    \ce{C} & complex \\
    \ce{E} & enzyme    \\
    $k_r$ & reaction rate coefficient \\
    $\MM$ & the set of labels of species  \\
    $\MM_0$ & the set of labels of initial species  \\
    $m$ & label of species  \\
   $M$ & number of species  \\
   $\N$& the set of positive integers\\
    \ce{P} & product\\
    $\RR$ & the set of labels of reaction steps \\
    $r$ & label of reaction step  \\
    $R$ & number of reaction steps\\
    $\R$& the set real numbers\\
    $\Rp$& the set positive real numbers\\
    \ce{S} & substrate \\
    \ce{X} & species \\
    \ce{X($m$)} & species  \\
    $y(t)$ & concentration of species \ce{Y} at time $t$  \\
    \ce{Y} & species  \\
    $\alpha(m,r)$ & reactant stoichiometric coefficient  \\
    $\beta(m,r)$ & product stoichiometric coefficient \\
    $\mu$&the largest \Vi\ of species\\
    $\varrho$&the largest \Vi\ of reaction steps\\
    \hline
  \end{tabular}
\end{table}
\section{Problem Formulation. Volpert Indices}\label{sec:problemform}
Here we shortly introduce the used formalism, the details can be found e.g.\ in \cite{feinbergbook} or \cite{tothnagypapp}.
Suppose we have a reaction vessel of constant volume, pressure and temperature and suppose it contains the chemical species $\ce{X}(1),\ce{X}(2),\ldots\ce{X}(M)$, where $M\in\N$ is a positive integer, and suppose that among the species the following reaction  steps take place:
\begin{equation}\label{eq:ccr}
  \sum_{m\in\MM}\alpha(m,r)\ce{X}(m)\longrightarrow\sum_{m\in\MM}\beta(m,r)\ce{X}(m)
  \quad(r\in\RR)
\end{equation}
with $\MM:=\{1,2,\cdots M\}$ and $\RR:=\{1,2,\cdots R\},$ where $R\in\N;$
and the elements of the matrices $\alphab$ and $\betab$ are supposed to be nonnegative integers, they are called \emp{stoichiometric coefficients}.
The object of our investigation is the mass action type kinetic differential equation of the reaction
\begin{equation}\label{eq:ikde}
\dot{c}_m=\sum_{r\in\RR}(\beta(m,r)-\alpha(m,r))k_r\prod_{p\in\MM}c_p^{\alpha(p,r)}
\quad(m\in\MM)
\end{equation}
describing the time evolution of the vector of \cc s $\cb=[c_1,c_2,\cdots c_M]\T$ of the species, where
$k_1,k_2,\cdots,k_R\in\Rp$ are the \emp{reaction rate coefficient}s.

\begin{figure}[!ht]
  \centering
  \includegraphics[width=0.8\textwidth]{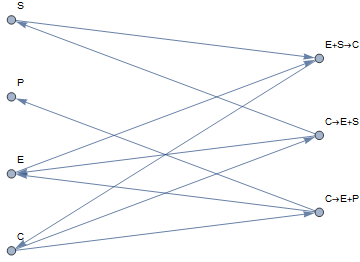}
  \caption{The \Vg\ of the Michaelis--Menten reaction}\label{fig:mmvg}
\end{figure}

There are different graphs representing the reaction
\cite[Chapter 3]{tothnagypapp}
of which we shall use the one actually introduced by Petri in 1939.
However, it was A. I. Volpert \cite{volpert} who started to connect properties of this graph to the properties of the solutions of induced kinetic differential equations. Further applications of this graph on the qualitative behavior of the solutions of \eqref{eq:ikde} can be found in \cite{banajicraciunmultiple,craciunpanteasontag,donnellbanaji,feliuwiuf,gaborhangosszederkenyi}, or on the decomposition of overall reactions \cite{kovacsvizvaririedeltoth}. (One can also meet rediscoveries in the literature: \cite{friedlertarjanhuangfan,reddymavrovouniotisliebman}.)
\begin{definition}
The Volpert graph of the reaction \eqref{eq:ccr} is a directed bipartite graph with the species and reaction steps as vertices, and with $\alpha(m,r)$ edges going from $\ce{X}(m)$ to $r$, and with $\beta(m,r)$ edges going from $r$ to $\ce{X}(m)$.
\end{definition}
Note that this is the graph most often used in papers and textbooks on me\-tab\-o\-lism. Fig. \ref{fig:mmvg} shows the \Vg\ of the Michaelis--Menten reaction.

Next we introduce indexing of the vertices of the \Vg.
\begin{definition}
Let $\MM_0\subset\MM$ be a subset of the species index set, called the indices of \emp{initial species}.
\begin{enumerate}
\item
Elements of $\MM_0$ receive zero index.
Zero will also be assigned to those reaction steps $r$ for which all the reactant species are available; in this case this means that all the reactant species are in $\MM_0$.
\item
Suppose that indices up to $k\in\N_0$ has been allotted. Then, a species receives the index $k+1$, if it has no index yet and there is a reaction step producing it.
A reaction step receives the index $k+1$ if it has no index yet and all the species needed for it to take place are available;
in this case this means that all the reactant species have an index not more than $k.$
\item
As the total number of vertices is finite, the algorithm ends after a finite number of steps. Species and reaction steps without index at termination will receive an index $+\infty$.
\end{enumerate}
\end{definition}
\begin{remark}
Notice that the value of the indices does not depend on the values of the stoichiometric coefficients, what only counts is their positivity.
I.e. if one has two reaction networks with stoichiometric coefficients that are positive for the same indices $m,r$ then the Volpert indices of the species and reaction steps will be the same. Our topics is quite different from the case when differences in the values of stoichiometric coefficients lead to qualitatively different consequences, e.g. \cite[p. 112]{feinbergbook}.
The time evolution of the concentrations also do depend on the values of the stoichiometric coefficients.
\end{remark}
Obviously, in the case when all the vertices receive a finite index, then the largest index $\varrho$ a reaction step can have is smaller than $R;$
the largest index $\mu$ a species step can have is
smaller than $M,$
and either $\varrho=\mu$ or $\varrho+1=\mu.$
The value of the \Vi\ shows a kind of distance of a vertex from the initial species, \cite[Theorem 9.1]{tothnagypapp}.
When applying the definition, the initial set of species will be those with positive \cc\ at the beginning.

\begin{figure}
  \centering
  \includegraphics[width=0.8\textwidth]{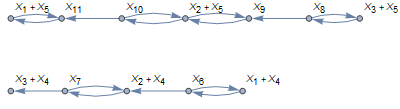}
  \caption{The \FHJg\ of the mitogen-activated protein kinase cascade}\label{fig:biofhjg}
\end{figure}

Let us see examples of indexing in the case of the Michaelis--Menten reaction.
\begin{example}
First, assume that \ce{E} and \ce{S} were initially present, then we get the following indices.

\begin{equation*}
\begin{tabular}{|l|c|}
  \hline
  Species & Indexes \\
  \hline
  \ce{S} & 0 \\
  \ce{E} & 0 \\
  \ce{C} & 1 \\
  \ce{P} & 2 \\
  \hline
\end{tabular}\quad
\begin{tabular}{|l|c|}
  \hline
  Reaction steps & Indexes \\
  \hline
  \ce{E + S -> C} & 0 \\
  \ce{C -> E + P} & 1 \\
  \ce{C -> E + S} & 1 \\
  \hline
\end{tabular}
\end{equation*}
Now a less favorable choice follows: the species initially present is \ce{P}.
\begin{equation*}
\begin{tabular}{|l|c|}
  \hline
  Species & Indexes \\
  \hline
  \ce{P} & 0 \\
  \ce{C} & $+\infty$ \\
  \ce{S} & $+\infty$ \\
  \ce{E} & $+\infty$ \\
  \hline
\end{tabular}\quad
\begin{tabular}{|l|c|}
  \hline
  Reaction steps & Indexes \\
  \hline
  \ce{C -> E + P} & $+\infty$ \\
  \ce{C -> E + S} & $+\infty$ \\
  \ce{E + S -> C} & $+\infty$ \\
  \hline
\end{tabular}
\end{equation*}
\end{example}
Let us consider a more complicated example.
\begin{example}
We take it from
\cite{markevichhoekkholodenko}, but it can also be found in the Biomodels Database;
see also \cite{lichtblau}.
\begin{equation*}
\begin{tabular}{|l|c|}
  \hline
  Species & Indexes \\
  \hline
  \ce{X_{10}} & 0 \\
  \ce{X_7}    & 0 \\
  \ce{X_{11}} & 1 \\
  \ce{X_5}    & 1 \\
  \ce{X_3}    & 1 \\
  \ce{X_2}    & 1 \\
  \ce{X_4}    & 1 \\
  \ce{X_9}    & 2 \\
  \ce{X_8}    & 2 \\
  \ce{X_1}    & 2 \\
  \ce{X_6}    & 3 \\
  \hline
\end{tabular}\quad
\begin{tabular}{|l|c|}
  \hline
  Reaction steps & Indexes \\
  \hline
  \ce{X_{10} -> X_{11}} & 0 \\
  \ce{X_{10} -> X2  + X5} & 0 \\
  \ce{X7  -> X3 + X4} & 0 \\
  \ce{X7 -> X2 + X4} & 0 \\

  \ce{X_{11} -> X1  + X5} & 1 \\
  \ce{X2 + X5  -> X_{10}} & 1 \\
  \ce{X2 + X5 -> X9} & 1 \\
  \ce{X3 + X5 -> X8} & 1 \\
  \ce{X2 + X4  -> X7} & 1 \\

  \ce{X1 + X5 -> X_{11}} & 2 \\
  \ce{X9 -> X2  + X5} & 2 \\
  \ce{X8  -> X9} & 2 \\
  \ce{X8 -> X3 + X5} & 2 \\
  \ce{X1 -> X4  + X6} & 2 \\

  \ce{X6  -> X2 + X4} & 3 \\
  \ce{X6  -> X1 + X4} & 3 \\
  \hline
\end{tabular}
\end{equation*}
where the meaning of the species are as follows.
MAPK is the abbreviation for mitogen-activated kinase in the table.
\begin{equation*}
\begin{tabular}{|l|l|l|}
  \hline
  Notation 1& Notation 2& Meaning \\
  \hline
  \ce{X_1}    & \ce{M}       &\text{dephosphorylated MAPK}\\
  \ce{X_2}    & \ce{Mp}      &\text{MAPK phosphorylated}\\
  &&\text{on one residue}\\
  \ce{X_3}    & \ce{Mpp}     &\text{MAPK phosphorylated}\\
  &&\text{on both residues}\\
  \ce{X_4}    & \ce{MAPKK}   &\text{MAPK kinase}\\
  \ce{X_5}    & \ce{MKP}     &\text{MAP kinase phosphatase}\\
  \ce{X_6}    & \ce{M-MAPKK} &\text{Complex}\\
  \ce{X_7}    & \ce{Mp-MAPKK}&\text{Complex}\\
  \ce{X_8}    & \ce{Mpp-MKP} &\text{Complex}\\
  \ce{X_9}    & \ce{Mp-MKP}  &\text{Complex}\\
  \ce{X_{10}} & \ce{Mp-MKP*} &\text{Complex}\\
  \ce{X_{11}} & \ce{M-MKP}   &\text{Complex}\\
  \hline
\end{tabular}
\end{equation*}
\end{example}
Now we give another kind of representation of this reaction, this time by the \FHJg\
shown in Fig. \ref{fig:biofhjg}. This representation can surely be reproduced by the reader without a formal definition; one should only remember that all the \emp{complexes} (formal linear combinations on both sides of the reaction arrows) should only appear once in the figure. Although all our calculations rely on the Volpert graph, the \FHJg\
may be more transparent in many cases.
In order to formulate our problem we need a theorem by Volpert, see e.g.\  \cite[Theorems 8.13 and 8.14]{tothnagypapp}.
\begin{theorem}
Suppose a set of species $\MM_0\subset\MM$ with positive initial \cc s is given.
\begin{enumerate}
\item
For all species $\ce{X}(m)\in\MM$ with a finite index and  for all times $t$ in the domain of existence of the solution of the \ikde\ $c_m(t) > 0$ holds.
\item
If the zero complex is not a reactant complex in the reaction \eqref{eq:ccr}, then for all species $\ce{X}(m)\in\MM$ with an infinite index and for all times $t$ in the domain of existence of the solution $c_m(t) = 0$.
\end{enumerate}
\end{theorem}
Cf. the content of the theorem to the behavior of absolute probabilities of a Markovian pure jump process (\cite[Remark 2.48]{liggett}), but  note that it is a slightly different  problem form the mathematical point of view: it relates an infinite system of linear differential equations while \eqref{eq:ikde} is a finite system of nonlinear equations.

If the zero complex is a reactant complex then the function \mma{VolpertIndexing} treats it as if it were a species automatically present in the set of initials, and the first part of the theorem still applies.
The result of
$$
\mma{VolpertIndexing[\{0 -> X\}, \{{ }\}, Verbose -> True]}
$$
\noindent is as follows.
\begin{equation*}
\begin{tabular}{|c|c|}
  \hline
  Species & Indexes \\
  \hline
  \ce{0} & 0 \\
  \ce{X} & 1 \\
  \hline
\end{tabular}\quad
\begin{tabular}{|c|c|}
  \hline
  Reaction steps & Indexes \\
  \hline
  \ce{0 -> X} & 0 \\
  \hline
\end{tabular}
\end{equation*}
showing that the (single) reaction step receives the index 0, because it starts from 0, supposed to be present.

Now our problem can be reformulated as finding such minimal initial set(s) with which the indices of all the species will be finite.
\section{Brute Force Tamed}\label{sec:bruteforce}
First, we describe and apply the brute force method, then we analyze its performance bounds.
\subsection{Steps of the brute force method}
Our first strategy is to form all possible nonempty subsets of $\MM$ (the power set except the empty set), and select those with which as initial species all the species have a finite \Vi.
Then select the minimal sets from the above set.
Although our calculations rely on the Volpert graph, the \FHJg\
is more transparent in many cases.
Selecting the minimal set means that given a set of sets we construct the underlying Sperner family (antichain or clutter) from it.
(Sperner family is a set of subsets of a finite set such that none of the element of the family contains another element of the family as a subset \cite{lovasz}.)
An easy to read Wolfram Language code for this (\mma{MinimalOnes}) can be found in the Appendix, which works as expected.
\begin{align*}
&\mma{MinimalOnes[\{\{C\}, \{E, S\}, \{E, C\}, \{S, C\}, \{C, P\},}\\
&\mma{\{\{E, S, C\}, \{E, S, P\}, \{E, C, P\}, \{S, C, P\},}\\
&\mma{\{E, S, C, P\}\}]}
\end{align*}
gives \mmares{\{\{C\},\{E,S\}\}}.
Now suppose we are given a reaction network and a set of initial species $\MM_0$ and we wonder whether all the species are generated from $\MM_0.$
This happens if all the \Vi es of the species are finite.
$$
\mma{FiniteVolpertIndexQ[\{"E","S"\}, "Michaelis-Menten"]}
$$
gives the expected result: \mmares{True}, whereas
$$
\mma{FiniteVolpertIndexQ[\{"E"\}, "Michaelis-Menten"]}
$$
will obviously be \mmares{False}.

Before calculating the power set of the species we introduce a few selection principles.
\begin{enumerate}
\item
If the reaction is not reversible then product species which are not reactant species at the same time, need not be taken into consideration. An example is the species \ce{P} in the Michelis--Menten reaction.
This is a requirement that can be realized in a purely formal, algorithmic way; from now on we shall start with the reactant species only. These are determined by \mma{ReactantSpecies}.

The result of \mma{ReactantSpecies["Michaelis-Menten"]} will be \mmares{\{C,E,S\}} as expected.
\item
The help of the chemist is needed to eliminate \emp{intermediates}, i.e.\ species which cannot be present initially because of chemical reasons. This set should always be specified by the user, although we could include such species as \ce{OH} forever, or we want to include \ce{C} in case of the Michaelis--Menten reaction. As to the second case: it would be desirable to algorithmically find such kind of intermediates, because they are defined by the structure of the reaction and not by the chemical nature of \ce{C}.
\item
A formal, automatic treatment can also be used to throw away \emp{subsets} of species where not all the atoms present in the reaction are contained in the species of the set.

A working code can again be found in the Appendix. If
\begin{align*}
&\mma{EmanuelKnorre = }\\
&\mma{\{Cl$_2\to$2Cl$^*$,CH$_4$+Cl$^*\to$ $^*$CH$_3$+HCl,$^*$CH$_3$+Cl$_2\to$CH$_3$Cl+Cl$^*$\}}
\end{align*}
then
\mma{AtomsPresentQ[EmanuelKnorre, \{HCl\}]} will give \mmares{False}, while
$$
\mma{AtomsPresentQ[EmanuelKnorre, \{HCl, C\}]}
$$
provides \mmares{True}.
\end{enumerate}
Now one can collect the species which can be present in the initial sets.
\begin{align*}
&\mma{PowersetBase[reac\_,intermediates\_:\{\}] :=}\\
&\mma{Complement[ReactantSpecies[\{reac\}], intermediates]}
\end{align*}
The result of \mma{PowersetBase["Michaelis-Menten", \{"C"\}]}
is \mmares{\{"E","S"\}} as expected, because \ce{P} as a terminal species has been discarded by the program, and the chemist asked us to discard the intermediate \ce{C}.
When forming the subsets of species one applies the check for the presence of all the necessary atoms.
The final code has been collected in the Appendix.

Here we give a few results, our notebooks contain a large set of examples.
$$
\mma{MinimalInitials["Michaelis-Menten"]}
$$
gives \mmares{\{\{"C"\},\{"E", "S"\}\}}, and $$\mma{MinimalInitials["Michaelis-Menten",\{"C"\}]}$$
gives \mmares{\{\{"E", "S"\}\}}.

\noindent\mma{MinimalInitials[EmanuelKnorre]},{\ }
\mma{MinimalInitials[EmanuelKnorre,\{\}]},\\
finally \mma{MinimalInitials[EmanuelKnorre,\{\},Atomic$\to$True]} all give
\mmares{\{\{CH$_4$,Cl$_2$\}\}}, but evaluation of the last command takes (in this case: only slightly) less time.

Finally, a larger example follows.

\noindent\mma{MinimalInitials[MAPKBIOMD26]} provides multiple solutions:
\begin{eqnarray*}
&
\{\ce{X_{10}},\ce{X_4}\},
\{\ce{X_{10}},\ce{X_6}\},
\{\ce{X_{10}},\ce{X_7}\},
\{\ce{X_{11}},\ce{X_4}\},
\{\ce{X_{11}},\ce{X_6}\},
\{\ce{X_{11}},\ce{X_7}\},&\\
&
\{\ce{X_4},\ce{X_8}\},
\{\ce{X_4},\ce{X_9}\},
\{\ce{X_5},\ce{X_6}\},
\{\ce{X_5},\ce{X_7}\},
\{\ce{X_6},\ce{X_8}\},
\{\ce{X_6},\ce{X_9}\}&\\
&
\{\ce{X_7},\ce{X_8}\},
\{\ce{X_7},\ce{X_9}\},
\{\ce{X_1},\ce{X_4},\ce{X_5}\},
\{\ce{X_2},\ce{X_4},\ce{X_5}\},
\{\ce{X_3},\ce{X_4},\ce{X_5}\}.&
\end{eqnarray*}

Further examples can be found in the notebook containing all the calculations.

\subsection{Performance bounds of the brute force method}
Let us start with the hardware. The test problems were solved on a laptop with 16 GB RAM and
Intel\textregistered\  Core\textsuperscript{TM} i7-8550U CPU \@ 1.80 GHz 1.99 GHz running a 64-bit Windows 10 operating system.

Two reasons may limit the method: the time needed to do the calculations and the memory needed. It turns out immediately that the second one is the real problem, therefore we dwell on this longer.

\begin{figure}[!ht]
  \centering
  \includegraphics[width=0.8\textwidth]{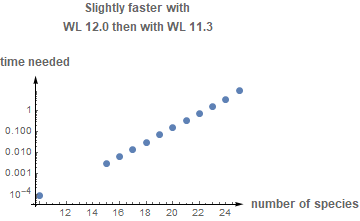}
  \caption{Times needed (from measurements) to form the power sets of sets having different sizes}\label{fig:times}
\end{figure}

First, we measured the time needed to form the power set of sets having different sizes: Fig. \ref{fig:times} It will turn out later that the reason we could not measure a timing for the 26 element set is that there is not enough memory space to store its power set.

\begin{figure}[!ht]
  \centering
  \includegraphics[width=0.8\textwidth]{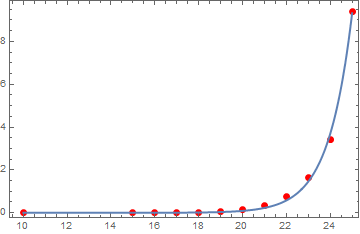}
  \caption{Curve fitted to measured times needed to form the power sets of sets having different sizes}\label{fig:fitting}
\end{figure}
(We have also found that Wolfram Language 12.0 is slightly faster than version 11.3.) Extrapolation from fitting a curve to the data (Fig. \ref{fig:fitting}) shows that with 30 element to form the subsets 16 minutes would be needed, with 40 elements 4 months, and with 50 elements 3500 years, a bit too long time.

The typical value of the variable \mma{MemoryAvailable[]} is \mmares{12,747,689,984} on our computer, and
\mma{Total[Map[ByteCount, Subsets[Range[25]]]]]} is
\mmares{11,506,401,520},
thus, without other ideas reactions with more than 25 species cannot be treated in such a simple way.

The final consequence should be that it is worth constructing other algorithms, as well.
\section{Using Integer Linear Programming}\label{sec:lp}
First, we have to introduce a series of further notations. 
\subsection{Notations}
Let $\Reactants(r)$ denote the species in the reactant complex of reaction step $r\in\RR:$
$$
\Reactants(r):=\{\ce{X}(m)\in\MM;\alpha(m,r)>0\},
$$
and $\Products(r)$ the species in the product complex of reaction step $r\in\RR:$
$$
\Products(r):=\{\ce{X}(m)\in\MM;\beta(m,r)>0\}.
$$
Furthermore, let $\Producing(m)$ be the set of those reaction steps which produce the species $m\in\MM:$
$$
\Producing(m):=\{r\in\RR;\beta(m,r)>0.\}
$$
We have a few binary variables (variables which can only have the value 0 or 1):
$y(m,\ell)=1$ if and only if species $\ce{X}(m)$ has a \Vi\ $\ell$, otherwise it is zero; $z(r,\ell)=1$ if and only if reaction step $r\in\RR$ has a \Vi\ $\ell-1$.
\subsection{Formulation of the problem}
Here we sketch how our problem can be attacked using the methods of integer programming problems. After the presentation of the method we analyse what is missing and what is different in this formulation compared to the original problem.

Find the minimum of the number of initial species
\begin{equation*}
{\mathcal O}(\yb,\zb):=\sum_{m = 1}^My(m, 0)
\end{equation*}
in the variables $\yb,\zb$ under the conditions \eqref{eq:cond1}--\eqref{eq:cond9} below. The meaning of the value of the objective function ${\mathcal O}$ is the total number of different initial species.
\begin{align}
&\forall m\in\MM \forall\ell\in\RR:y(m,\ell)\mbox{ is }0\mbox{ or }1;
\label{eq:cond1}\\
&\forall m\in\MM \forall \ell\in\RR:z(r,\ell)\mbox{ is }0\mbox{ or }1.
\label{eq:cond2}
\end{align}
Eqs. \eqref{eq:cond1} and \eqref{eq:cond2} express the requirement that the variables are binary variables.
\begin{equation}
\forall m\in\MM :\sum_{\ell=0}^{R-1}y(m,\ell)=1;
\label{eq:cond3}\end{equation}
Condition \eqref{eq:cond3} means that all the species are produced with some finite Volpert index.
\begin{equation}
\forall r\in\RR :\sum_{\ell=0}^{R-1}z(r,\ell)=1;
\label{eq:cond4}\end{equation}
Eq. \eqref{eq:cond4} means that all the reaction steps take place (i.e.\ they have a finite Volpert index). Note that if the \Vi\ of all the species is finite then no reaction step can have an infinite \Vi.

From now on the absolute value sign is used to denote the number of elements of a finite set.
\begin{equation}
\forall r\in\RR \forall k\in\RR :\sum_{m\in \Reactants(r)}\sum_{\ell=0}^ky(m,\ell)-\sum_{\ell=1}^{k}z(r,\ell)
\le |\Reactants(r)|-1;
\label{eq:cond5}
\end{equation}
If all the reactants are present, then the reaction step can occur.

\begin{equation}
\forall r\in\RR \forall k\in\RR :\sum_{m\in \Reactants(r)}\sum_{\ell=0}^{k-1}y(m,\ell)-|\Reactants(r)|\sum_{\ell=1}^{k}z(r,\ell)
\ge0;
\label{eq:cond6}
\end{equation}
If any of the reactant is absent, then the reaction step cannot occur.
\begin{equation}
\forall r\in\RR \forall k\in\RR \forall m\in \Products(r):-\sum_{\ell=0}^{k}y(m,\ell)+z(r,k)\le 0;
\label{eq:cond7}
\end{equation}
If the reaction step $r$ has the \Vi\ $k,$ then at least one of the species has a \Vi\ not larger than $k.$
\begin{equation}
\forall r\in\RR:\sum_{\ell=1}^{R}z(r,\ell)\ge1 ;
\label{eq:cond8}
\end{equation}
The \Vi\ of all the reaction steps is finite.

\begin{equation}\forall m\in\MM \forall k\in\RR:\sum_{\ell=1}^ky(m,\ell)\le\sum_{r\in\Producing(m)}\sum_{\ell=1}^kz(r,\ell).
\label{eq:cond9}\end{equation}
The Volpert indices of the product species of a reaction step are not greater than the \Vi\ of the reaction step plus 1.

First, note that the minimum of the number of initial species does not define the same initial set that we called minimal above. Second, further conditions can be formulated to eliminate end products and intermediates, etc.
Instead of refining this approach we provide another one that works for large data sets and fits most naturally to the structure of the problem.
\section{The Smart Solution}\label{sec:smart}
First, we describe lexicographic and reverse lexicographic orderings, then we show how to apply them to the solution of our problem.
\subsection{The lexicographic and reverse lexicographic orderings as tools of enumeration}
To obtain the solution of many different problems
a large number of objects must be enumerated.
The enumeration must be fast and exact.
An enumeration is fast if no object is considered more
than once and it is exact if all relevant objects are considered.
Brute force is not fast as it considers many irrelevant objects.
Thus, the enumeration must consider every investigated object exactly once.
Notice that this kind of enumeration gives a complete order of the objects.
In a complete order, any two different objects are comparable, i.e. one of them precedes the other one.

There are two types of constituents in a chemical reaction network, species and reaction steps.
(We may---and shall---also take into consideration that species are built up from atoms.)
The problem of determining the minimal sets of initial species in such a way that starting from these sets, the reaction steps produce all species, needs the enumeration only one of the two subsets, i.e. either the subsets of species, or the subsets of reaction steps.
We shall be working with the first.

These subsets, i.e. the subsets of species are the objects to be enumerated.
Assuming that there are $M:=|\MM|$ species,
the total number of all the subsets is $2^M$, which can be very huge.
Taming brute force meant to exclude some species or some subsets of species based on preliminary considerations, not an essential speed-up.
Thus, an effective enumeration method must handle the majority of cases in an implicit way.
It means that the individual investigation of the implicitly enumerated cases is not necessary.
For example, if a subset $T$ satisfies all the requirements,
then the investigation of any subset $U$ such that $T\subset U$ is superfluous.

The subsets must be represented in a way which is suitable for calculation.
The \emp{characteristic vector} of the subsets is one possible tool.
The characteristic vector $u$ of a subset $\MM_0\subset\MM$ is an
$M$-dimensional vector such that
\begin{eqnarray*}
u_m:=\left\{
\begin{array}{ll}1, & \text{ if } m\in \MM_0; \\ 0,&\text{ if } m\not\in \MM_0.\end{array}
\right.
\end{eqnarray*}
Such vectors can be ordered by the lexicographic order of vectors.
A vector such that every component is either 0, or 1, is called
\textbf{binary} or \textbf{0-1} or \textbf{Boolean vector}.
The number of binary vectors of fixed dimension is finite.
Let us introduce a lexicographically increasing order of such vectors.
\subsubsection{Lexicographic ordering.}
\begin{definition}
Let $u,w \in \{0,1\}^M$ be two $M$-dimensional vectors.
Vector $u$ is \emp{lexicographically smaller} than vector $w$ if an index $p: 1< p\leq M$ exists such that
\begin{enumerate}
    \item  $u_m=w_m  \mbox{ for }m=1,\ldots,p-1$, and
    \item $u_p<w_p$.
\end{enumerate}
This relation is denoted by $u\prec w$.
\end{definition}
Thus, lexicographically the smallest vector is $(0,0,\dots,0,0)$, and the second smallest vector is $(0,0,\dots,0,1)$.
Assume that $u$ is a non-zero binary vector, different from the one only having 1 as coordinates.
Let $p$ ($1\le p < M$) be the last index such that $u_p=0, u_{p+1}=1.$
In the lexicographical order, the \emp{next} vector is $v$ with
$u_m=v_m$  $(m=1,\dots,p-1),v_p=1,v_{p+1}=\dots v_M=0.$
Consequently, the vector such that all coordinates are 1, is the lexicographically largest binary vector, and no vector exists thereafter.

Note that lexicographic ordering orders the binary vectors in the same way as if they were ordered according to the values of the number composed from the components as digits in the binary system.
\subsubsection{Reverse lexicographic ordering.}
Another ordering, called \emp{reverse lexicographic ordering} is obtained if the coordinates are used in the opposite order.
\begin{definition}
Let $u,w \in \{0,1\}^M$ be two $M$-dimensional vectors.
Vector $u$ is \emp{smaller with respect to the reverse lexicographic ordering} than vector $w$ if
an index $p: 1< p\leq M$ exists such that
\begin{enumerate}
    \item  $u_m=w_m  \mbox{ for }m=p+1,\ldots,M$, and
    \item $u_p<w_p$.
\end{enumerate}
\end{definition}
If the reverse lexicographic order of the coordinates
is used, then the rule is as follows:
the first 0 from left is changed to 1 and all coordinates before it must be zero.
Assume that the current vector
is $u$ with $u_1=\dots=u_p=1, u_{p+1}=0$.
Subsequently, the next vector $v$ is such that $v_1=\dots=v_p=0, v_{p+1}=1$ and for $m=p+2,...,M: v_m=u_m.$
Specifically, $(0,0,...,0,0)$ is the smallest, and $(1,1,...,1,1)$ is the largest vector
also in the case if reverse lexicographic order is applied.
However, the second and third smallest
vectors are $(1,0,...,0,0)$, and $(0,1,...,0,0)$, respectively.
In what follows, this \emp{reverse lexicographic order} will also be used.
See Table \ref{tab:enumeration}.
\begin{table}[!ht]
\begin{center}
\begin{tabular}{cccc}
Case &	Vector &	Case &	Vector\\
1 &	(0,0,0,0) &	9 &	(0,0,0,1)\\
2 &	(1,0,0,0) &	10 &	(1,0,0,1)\\
3 &	(0,1,0,0) &	11 &	(0,1,0,1)\\
4 &	(1,1,0,0) &	12 &	(1,1,0,1)\\
5 &	(0,0,1,0) &	13 &	(0,0,1,1)\\
6 &	(1,0,1,0) &	14 &	(1,0,1,1)\\
7 &	(0,1,1,0) &	15 &	(0,1,1,1)\\
8 &	(1,1,1,0) &	16 &	(1,1,1,1)\\
\end{tabular}
\caption{The enumeration of the 16 binary vectors of the 4-dimensional space if the reverse lexicographic order is used.}\label{tab:enumeration}
\end{center}
\end{table}
Using again the Michaelis--Menten example we find that Case 5 precedes Case 6, meaning that the potential initial set $\{\ce{C}\}$ precedes the initial set $\{\ce{E},\ce{C}\}.$

Note that beyond the two orderings above any other order of the coordinates gives another lexicographic order.
\subsection{The application of reverse lexicographic ordering to the enumeration of minimal sets of species}
The logic of the enumeration in the investigation of a chemical reaction is as follows.
Two actions must be done with every vector \(u\).

The first one is the checking if all materials have a finite index in the Volpert graph given that $u$ is the characteristic vector of initial species.
One possible answer to this question is given
by  the procedure \emp{complete} in the pseudocode of Fig. \ref{fig:lexicographic},
or the above used Boolean function
\mma{FiniteVolpertIndexQ} defined in the Appendix \ref{sec:appendix}.
Its value is true if the condition is satisfied, otherwise the value is false.

The second one is that the vector immediately following $u$ is to be generated. It is done by the function $\following{u}$.
The minimal solutions are collected in set $U.$
The algorithm puts a vector $u$ into this set, if the vector represents a subset as its characteristic vector which is currently minimal. If there are vectors in $U$ which represent subsets larger than the subset represented by $u$, then they must be deleted from $U.$

The rule how to obtain the lexicographically next vector is as follows:
Find the last 0 in the vector, i.e. the 0 component with the highest index.
Thus, the end of the vector is $(\dots,0,1,\dots,1)$.
Change the zero to 1 and behind it all components become 0.
The code for \following\ { }in the case of reverse lexicographic ordering is similar.
\begin{figure}[!ht]
\textbf{1 for $m:=M$ step $-1$ to $1$ do\\
2\quad  if $u_m=0$\\
3\quad  then\\
4\quad  begin\\
5\quad\quad      $u_m:=1$\\
6\quad\quad for $p:=m+1$ to $M$ do $u_p:=0$\\
7\quad end}
\caption{This is the code for the function \emp{next} with lexicographic ordering}
\label{fig:next}
\end{figure}

\begin{figure}
\textbf{{1} begin\\
{ 2} \quad 	$u:=(0,0,\dots,0,0)$\\
{ 3} \quad $U:=\emptyset$\\
{ 4} \quad 	while $u\prec (1,1,\dots,1,1)$ do\\
{ 5} \quad  begin\\
{ 6} \quad\quad  if complete($u$)\\
{ 7} \quad\quad  then if $V\,:=\,\{v\in U:u\prec v\}\neq\emptyset$\\
{ 8} \quad\quad\quad\quad  	then $U\,:=\, U\cup\{u\}\setminus V$\\
{ 9} \quad\quad\quad\quad  	else if \quad$\not\exists v\in U:v\prec u$\\
10 \quad\quad\quad\quad\quad\quad  then $U\,:=\,U\cup\{u\}$\\
11 \quad\quad  $u\,:=\,\following{u}$\\
12 \quad  end\\
13 end}
    \caption{The frame of the lexicographic enumeration}
    \label{fig:lexicographic}
\end{figure}
With the code in Fig. \ref{fig:lexicographic} the problem has been solved, unless one wants to take into consideration the atomic structure of species to reduce the time needed to the calculations.
\subsection{Algorithmic details of the lexicographic enumeration in presence of atoms}
\subsubsection{The order of coordinates applied in this description.}
How to enlarge the algorithm of Fig. \ref{fig:lexicographic} if atomic structure of the species is also considered?

First, the lexicographic order is applied.
It means that the lexicographically smallest vector is $(0,0,\dots,0,0,0)$, the second smallest is $(0,0,\dots,0,0,1)$, the third smallest is $(0,0,\dots,0,1,0)$, etc.

The dimension of the vector, i.e. $M$, is the number of species in our case. Index $m$ identifies the species \ce{X($m$)}.
The order of the coordinates is important from algorithmic point of view. 
\subsubsection{Checking the atomic condition (A).}
Let $\AAA$ be the set of the different atoms in the chemical species. Similarly, let $\MM$ be the set of species.
The presence of atom $a$ is described by a characteristic vector as follows.
The vector
\[c^a\in \{0,1\}^M\]
is the characteristic vector of the atom $a\in\AAA$ if its components are defined as follows:
\begin{eqnarray*}
c_m^a=\left\{\begin{array}{ll}1,&\hbox{if atom } a \hbox{ is present in species } \ce{X($m$)}\\ 0&\hbox{otherwise.}\end{array}\right.
\end{eqnarray*}
These vectors are fixed during the whole process. Let $L^a$ be defined as \[L^a:=\max\{m|c_m^a=1\},\]
i.e. $L^a$ is the last species to contain atom $a.$
Let $u$ be the characteristic vector of the initial set under construction, i.e. if $u_m=1$, then species $\ce{X($m$)}$ is present in the set, otherwise 0.
An atom is represented in the set, if there is a species in the set
$\{\ce{X(m)};u_m=1\}$ such that it contains atom $a$, or $c_m^a=1.$
Hence, an atom is represented in the set if and only if
\[\sum_{m=1}^M c_m^a u_m \geq 1.\]
If this inequality is not satisfied, then the equation
\[\sum_{m=1}^M c_m^a u_m=0\]
holds implying $u_{L^a}=0$.
Then, species $\ce{X($L^a$)}$ is put into the initial set, i.e. $u_{L^a}$ is changed to 1.
The result is the lexicographically minimal set which contains all the species being present at the beginning of the construction of the current suggestion for the initial set and contains atom $a$ as well.
These steps are done for all kinds of atoms.
\subsubsection{Checking the minimality condition (B).}
Assume that the set of the characteristic vectors of the minimal sets is set $U$.
A smaller than $u$ initial set is contained in $U$ if and only if there is a vector $v\in U$ such that $v\not\prec u.$
Notice, that the two vectors cannot be equal as vector $u$ is lexicographically strictly increasing during the process.
In other words, this condition means that $v_m=1$ implies  $u_m=1$, or equivalently
\[\sum_{m=1}^M v_m u_m = \sum_{m=1}^M v_m.\]
If the condition is satisfied, then the current candidate cannot be minimal.
Assume that $\mathrm{max1}$ is the index of the last component having value 1 in the current candidate vector $u$. Further on, let $\mathrm{max0}$ be the index of the last 0 component in $u$ before the $\mathrm{max1}$ component. Thus
\[\mathrm{max1}=\max\{m|u_m=1\} \hbox{ and } \mathrm{max0}=\max\{m|m<\mathrm{max1}, u_m=0\}.\]

Thus, the end of the candidate vector is 01\dots10\dots0. If any 0 component of the last segment of the vector consisting of only 0s is changed from 0 to 1, then the non-minimality condition still holds. Thus, component $\mathrm{max0}$ must be changed from 0 to 1 as the lexicographically minimal attempt to get out from the non-minimality state. All components behind this component must be 0.
It is one of the ways how lexicographic order can enumerate many cases
in an implicit way.

\subsubsection{The order of the components.}
Notice, that $u_m=1$ in every second iteration in the case of a complete enumeration. Therefore, it is advantageous if it represents a species which is present in a minimal set with high probability.
On the other hand, $u_1=0$ in the first half of the complete enumeration.
Thus, if any minimal set exists containing the first species, it can be obtained only in the second half of the enumeration. To have many minimal sets can accelerate the enumeration discussed in the previous sections.
Thus, if there is any available information on the probability of the presence of the species in minimal initial sets, it is worth to change the order of the species accordingly.

\subsubsection{Algorithmic summary.}
Here is a short algorithmic summary of the discussed procedure. The Boolean variable Volpert-able is true if the current candidate vector satisfies the necessary conditions (A), and (B) and the construction of the Volpert-graph may begin. No new notation is introduced. Relation = means logical condition while  :=  means assignment, as usual.
\begin{figure}[!ht]
\textbf{{1} Volpert-able:=false\\
2 while not(Volpert-able) do\\
3 begin\\
4 \quad  for $a\in A$ do\\
5 \quad\quad if $\sum_{m=1}^n c_m^a u_m=0$ then $u_{L^a}:=1$\\
6 \quad Volpert-able:=true\\
7 \quad for $v\in U$ do\\
8 \quad\quad if $v\not\prec u$ then\\
9  \quad\quad begin\\
10 \quad\quad\quad  Volpert-able:=false\\
11 \quad\quad\quad  $\mathrm{max1}:=\{m\mid u_m=1\}$\\
12 \quad\quad\quad  $\mathrm{max0}>=\{mm<\mathrm{max1},u_m=0\}$\\
13 \quad\quad\quad  $u_{\mathrm{max0}}:=1$\\
14  \quad\quad\quad  for $m:=\mathrm{max0}$ to $\mathrm{max1}$ do $u_m:=0$\\
15  \quad\quad end\\
16 end}
\caption{Summary of the algorithm}
\end{figure}
\subsection{Results}
Several well known problems were solved using the proposed lexicographic enumeration method to analyse the efficiency of the algorithm.
The test problems were solved on a laptop with 8 GB RAM and Intel\textregistered\  Core\textsuperscript{TM} i7-8650U CPU \@ 1.90 GHz 2.11 GHz running a 64-bit Windows 10 operating system. The problems solved consisted of
hydrogen combustion models \cite{oconairecurransimmiepitzwestbrook}
and Leeds3 \cite{smithgoldenfrenklachmoriaryeiteneergoldenbergbowmanhansonsonggardinerlissianski},
methanol combustion models, such as
Klippenstein2011 \cite{klippensteinhardingdavistomlinskodje},
Li2007 \cite{lizhaokazakovchaosdryerscire} and
Rasmussen2008 \cite{rasmussenhansenmarshallglarborg}.
(See also \cite{tothnagyzsely} for a collection of combustion models and the references therein.)
The specifics of the problems such as total number of molecules, atoms and interactions alongside the total CPU times are reported in Table \ref{tab:smartresults}.

\begin{table}[!ht]
  \centering
\begin{tabular}{ |l|c|c|c|c|c| }
 \hline
 Ref. & Species & Reaction&Atoms& Inter- & CPU \\
 no.&&steps&&mediates&time\\
  \hline \hline
 \cite{zselyolmpalvolgyivarganagyturanyi} & 8 &40& 2& 4 & 0.012   \\
  \hline
 \cite{klippensteinhardingdavistomlinskodje} & 18 & 172&3 & 8 & 0.029   \\
 \hline
 \cite{smithgoldenfrenklachmoriaryeiteneergoldenbergbowmanhansonsonggardinerlissianski}&9&44&2&4&0.014 \\
 \hline
 \cite{lizhaokazakovchaosdryerscire}       &18&170& 3&  8&0.39 \\
 \hline
 \cite{rasmussenhansenmarshallglarborg}&28&320& 4& 13&26.361 \\
 \hline
\end{tabular}
\caption{Results of the smart method on a few large models} \label{tab:smartresults}
\end{table}

The algorithm fails to solve problems with large number of elements in an acceptable time. The reason can be easily attributed to the large set of possible binary vectors which amounts to $2^n$, where $n$ is the total number of species excluding the intermediate ones. As an instance, the problem in Ref. \cite{zabettahupa} contains 58 elements from which 37 elements are intermediates, which leaves the algorithms 21 elements to deal with. Thus, cardinality of the set of possible lexicographic vectors is $2^{21}=2,097,152$. Although this number is not that large, it takes more than 72 hours to solve this problem. One of the reasons is the time it takes for each vector to transform and exclude the intermediate elements. Consequently, one of the solutions lies in the first step of the solution method which is the coding of chemical species into numbers. If species other than intermediates are coded first, the lexicographic vector can be divided into two parts; first part describes the existence of main species and the second part which describes the intermediate ones. While calculating the next lexicographic order, the algorithm can easily dispose of the second part and all the related transformations which can effectively speed up the algorithm.
The second solution for large problems is limiting the search space to minimal sets with maximum cardinality of $n$. For properly small $n$, the cardinality of the set of possible lexicographic vectors will reduce to $c(m,n)$, where $m$ is the number of non-intermediate elements.
The results using the suggestions above is described in what follows.
By limiting the cardinality of minimal sets to 5, the problem of Zabetta and Hupa \cite{zabettahupa} was solved in just 780.49 seconds and 122 minimal sets were found.

The problem in Aranda et al. \cite{arandachristensenalzuetaglarborggersencgaodmarshall} with 1063 reactions and 76 molecules from which 40 were intermediates was solved within 10729.79 seconds. The maximum cardinality of minimal sets was set to 5 and 434 minimal sets were found.
\section{Conclusions and perspectives}

\subsection{Possible practical applications}
An obvious medical application of the solution of the problem is as follows. Suppose an organism should have all the species in a part of the metabolism, and one has to provide initial species to ensure this.
Then, we have to choose between the minimal initial species sets in such a way that we provide those species which cause the smallest side-effects, which are easiest to obtain, which are the cheapest, which has the shortest absorption time, etc.

In the case of chemical engineering applications (of which combustion is just one of the possible areas) one of the viewpoints can again be economics  or availability, but it may also happen that starting from one minimal set has smaller unwanted effect on the environment.
\subsection{Theoretical connections}
Let us start from an appropriate minimal initial set.
Then all the species will have positive concentrations for all finite (positive) times, or, the trajectory lies in the open first orthant of the species space. In a large class of reaction networks (for weakly reversible networks) this property (permanence) is conjectured to imply (\cite{brunnercraciun,craciunnazarovpantea}) that the same property holds for "infinite time", as well.
\subsection{Further problems}
\begin{enumerate}
\item
An extension for the case of reaction-diffusion models seems to be possible based on the papers by \cite{volpert} and the more recent ones by \cite{minchevaroussel,minchevasiegelstab,minchevasiegelpos}.
Exclusion of cross diffusion would not be a real restriction, if necessary from the mathematical point of view.
\item
In our Copenhagen lecture \cite{tothcopenhagen} we presented the following theorem.
\begin{theorem}\label{th:copenhagen}
Given a kinetic differential equation choose a few coordinates of the stationary point which is zero, put them into the equations. What you get will also be a kinetic differential equation.
\end{theorem}
Based on it one could select some coordinates to be zero, and construct appropriate inducing reactions to the remaining smaller kinetic differential equation, and apply one of the structural conditions to say something about the remaining coordinates if they are positive or zero.
\item
Beyond the Theorem \ref{th:copenhagen} above one might also study \emp{siphons}, starting with \cite{shiusturmfels} and continuing with recent developments.
\item
Which are the states which can be reached by a reaction at all?
This is connected to controllability, and also to the range of the right hand side etc.
\cite[p. 164]{tothnagypapp}.
\item
Beyond kinetics: given a polynomial differential equation we may be able to answer the analogous question, if we are able to find a reaction network inducing the give differential equation. Ref. \cite{craciunjohnstonszederkenyitonellotothyu} provides answers to questions of this last type.
\end{enumerate}

\section*{Conflicts of interest}
There are no conflicts to declare.

\section*{Acknowledgements}
The paper is based on a talk given by the authors at the 11th Conference on  Mathematical and Theoretical Biology, 23--27 July, 2018, Lisboa.
Discussions in the friendly and constructive atmosphere of the Workshop on the Advances in Chemical Reaction Network Theory
at the Erwin Schr\" odinger Institue, Vienna, significantly contributed to the paper.
JT thanks for the support of the National Research, Development and Innovation Office (SNN 125739).
Careful reading by Prof. Vilmos G\'asp\'ar and
Mr. Michael Gustavo were also a great help.

\bibliography{Positiv}
\bibliographystyle{plain}
\section{Appendix: The Code of the Brute Force Solution}\label{sec:appendix}
This code is for selecting the minimal subsets from the good ones.
\begin{align*}
&\mma{ProperSubsetQ[x\_, y\_] := }\\
    &\mma{And[SubsetQ[x, y], Not[Equal[x, y]]];}\\
&\mma{MinimalOnes[l\_List] := Fold[Function[\{x, y\},}\\
    &\mma{Select[x, Not[ProperSubsetQ[\#, y]] \&]], l, l]}
\end{align*}

Let us decide if $+\infty$ is among the \Vi\ in a reaction with a given set of initial species.
\begin{align*}
&\mma{FiniteVolpertIndexQ[ini\_, reac\_] :=}\\
&\mma{FreeQ[VolpertIndexing[{reac}, ini], Infinity, 3]}.
\end{align*}

It is only the reactant species which are to be considered.
\begin{align*}
&\mma{ReactantSpecies[reac\_] := Union[Cases[Level[Map[}\\
&\mma{Part[\#,1]\&, ReactionsData[\{reac\}]["fhjgraphedges"]],}\\
&\mma{\{-1\}],Except[\_Integer]]}.
\end{align*}

A possible set of intermediates may be
\begin{align*}
&\mma{inter=\{Br$^-$,BrO$_3^-$,CH$_3$,Cl$^*$,H,H$^+$,HO$_2$,I$^-$,IO$_3^-$,O,OH$^-$,}\\
&\mma{Mn(OH)$_2^+$,C\}}
\end{align*}
but the default is that there are no intermediates present.

The function \mma{GetElements} gives all the atoms present in a set of species (having atomic structure), and also, all the atoms present in a set of reaction steps. Both
$$
\mma{GetElements[\{H, Cl, HCl\}]}
$$
and
$$
\mma{GetElements[\{H + Cl -> HCl\}]}
$$ give \mmares{\{Cl, H\}}.
Using this function one can easily formulate the code needed.
\begin{align*}
&\mma{AtomsPresentQ[reac\_, subset\_] := }\\
&\mma{GetElements[reac] == GetElements[subset]}
\end{align*}

In case if the species have an atomic structure, we take this fact into consideration, otherwise not.
An option will be introduced to take care of this.
Next, we form the difference of the set of species and the intermediates.
In general, this set will be the base set from which to form the power set.
In the atomic case we restrict ourselves further to those sets which contain all the atoms present in the reaction network.
Then we select those elements of the power set which produce finite \Vi es, and finally select the minimal ones among them.

\begin{align*}
&\mma{Options[MinimalInitials] = \{Atomic $\to$ False\};}\\
&\mma{MinimalInitials[reac\_, intermediates\_: \{\},}\\
&\mma{OptionsPattern[]] := Module[\{MM, pb\},}\\
&\mma{pb = Complement[reactantSpecies[{reac}], intermediates];}\\
&\mma{MM = If[OptionValue[Atomic],}\\
&\mma{Select[Select[Rest[Subsets[pb]], }\\
&\mma{(GetElements[\#] == GetElements[pb])\&],}\\
&\mma{FiniteVolpertIndexQ[\#, reac]\&],}\\
&\mma{Select[Rest[Subsets[pb]],}\\
&\mma{(FiniteVolpertIndexQ[\#, reac])\&]]; MinimalOnes[MM]]}\\
\end{align*}

Dealing with the atoms first may strongly reduce the number of the elements which will form the power set.
In the case of the Briggs-Rauscher reaction it is only 37\% of the subsets that have to be considered, as the result of
\mma{AtomicSaving[GetReaction["Briggs-Rauscher"]]} is \mmares{0.372082}, using the code
\begin{align*}
&\mma{AtomicSaving[reac\_] := Module[}\\
&\mma{\{ini = Rest\@Subsets\@ReactionsData[\{reac\}]["species"], sel,}\\
&\mma{total = GetElements[\{reac\}]\},}\\
&\mma{sel = Select[ini, GetElements[\#] == total\&];}\\
&\mma{N[Length[sel]/Length[ini]]]}
\end{align*}

The ratio of subsets left in a few reaction can be seen in the Table \ref{tab:saving}.

\begin{table}[!ht]
  \centering
  \begin{tabular}{|l|l|}
     \hline
     Name           & Ratio of subsets left \\
     \hline
     Briggs-Rauscher& 0.37 \\
     Chapman cycle  & 0.47 \\
     Emanuel-Knorre & 0.84 \\
     FKN            & 0.26 \\
     Ogg            & 0.97 \\
     Petri          & 0.22 \\
     Vaiman         & 0.24 \\
     \hline
   \end{tabular}
  \caption{Saving via utilizing atomic composition}\label{tab:saving}
\end{table}

\end{document}